\documentclass[12pt]{amsart}
\usepackage{amssymb}
\usepackage{amsmath,amssymb,amsfonts,amsthm,graphics,
latexsym, amscd, amsfonts, epsfig, eepic,epic,psfrag,setspace}
\usepackage{mathrsfs}
\usepackage{color}
\usepackage{eucal}
\usepackage{eufrak}
\usepackage[all]{xypic}
\usepackage{xspace}

\theoremstyle{plain}
\newtheorem{theorem}{Theorem}

\newtheorem{lemma}{Lemma}

\theoremstyle{definition}

\theoremstyle{example}

\theoremstyle{remark}

\numberwithin{equation}{section}

\begin{document}
\doublespacing
\begin{center}
{\bf\Large Pseudoknot RNA structures with arc-length $\ge 4$}
\\
\vspace{15pt}Hillary S.~W. Han and Christian M. Reidys$^{\,\star}$
\end{center}

\begin{center}
Center for Combinatorics, LPMC-TJKLC\\
Nankai University  \\
Tianjin 300071\\
         P.R.~China\\
         Phone: *86-22-2350-6800\\
         Fax:   *86-22-2350-9272\\
reidys@nankai.edu.cn
\end{center}
\centerline{\bf Abstract}

\qquad \:In this paper we study $k$-noncrossing RNA structures with
minimum arc-length $4$ and at most $k-1$ mutually crossing bonds.
Let ${\sf T}_{k}^{[4]}(n)$ denote the number of $k$-noncrossing RNA
structures with arc-length $\ge 4$ over $n$ vertices. We prove (a) a
functional equation for the generating function $\sum_{n\ge 0}{\sf
T}_{k}^{[4]}(n)z^n$ and (b) derive for $k\le 9$ the asymptotic
formula ${\sf T}_{k}^{[4]}(n)\sim c_k\, n^{-((k-1)^2+(k-1)/2)}\,
\gamma_k^{-n}$. Furthermore we explicitly compute the exponential
growth rates $\gamma_k^{-1}$ and asymptotic formulas for $4\le k\le
9$.



\section{Introduction}\label{S:Introduction}

RNA pseudoknot structures \cite{Science:05a,Westhof:92a} are a
reality. They occur in functional RNA (RNAseP \cite{Loria:96a}),
ribosomal RNA \cite{Konings:95a} and are conserved in the catalytic
core of group I introns. Due to the crossings of arcs their theory
differs considerably from RNA secondary structures. Pseudoknots are
inherently noninductive and the standard dynamic programming folding 
paradigm employed for RNA secondary structures can only generate particular
subclasses of pseudoknot structures \cite{Rivas:99}. 
Recently the concept of $k$-noncrossing RNA structures has been 
introduced \cite{Reidys:07pseu}. Here the idea is that the complexity 
of the structure is captured by an inherently ``local'' property: the
maximal number of mutually crossing bonds. A structure is
$k$-noncrossing, if there exists no $k$-set of mutually crossings
arcs. The locality is in fact of central importance: point in case are
RNA bisecondary structures introduced by P.F.~Stadler \cite{Stadler:99}.
These structures are constructed as superpositions of two RNA secondary 
structures and correspond to {\it planar} $3$-noncrossing 
structures \cite{Reidys:07pseu}. The planarity property is clearly
non-local and at present 
time the generating function for RNA bisecondary structures is not known.

\begin{figure}[ht]
\epsfig{file=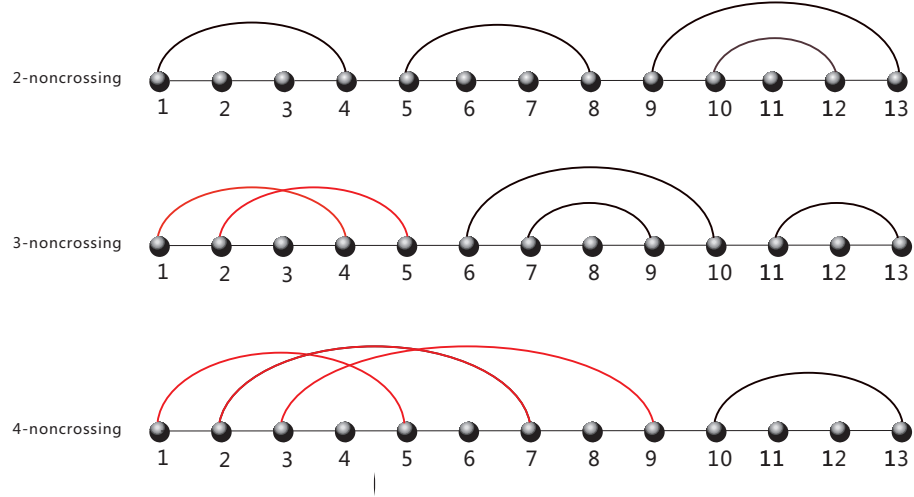,width=0.8\textwidth} \caption{\small
$k$-noncrossing structures: $2$- $3$- and $4$-noncrossing structures
(top to bottom). Maximal sets of mutually crossing arcs are colored
red.} \label{F:diagram}
\end{figure}

A very intuitive approach to the $k$-noncrossing property of RNA
molecules is their diagram representation \cite{Stadler:99}. It is 
obtained by drawing the nucleotide-labels $1,\dots,n$ in increasing 
order in a horizontal line and drawing the arc-labels $(i,j)$ in the 
upper half-plane, if and only if $i$ and $j$ are paired in the structure, 
see Figure~\ref{F:diagram}.
 We call a diagram $k$-noncrossing, if it
does not contain $k$ mutually crossing arcs. The length of an arc
$(i,j)$ is given by $\lambda=j-i$ and a stack of length $\sigma$ is
a sequence of ``parallel'' arcs of the form
$((i,j),(i+1,j-1),\dots,(i+(\sigma-1), j-(\sigma-1)))$.
A $k$-noncrossing RNA structure is a $k$-noncrossing diagram over $[n]$
having minimum arc-length $\lambda>1$. These structures have been studied in
\cite{Reidys:07pseu,Reidys:07asy1} via a bijection into
vaccillating tableaux in the context of tangled diagrams \cite{Reidys:07vac}.
For the enumeration of structures with crossing arcs the
tableaux-interpretation is non-optional. There is, to the best of our 
knowledge, no way to inductively construct $k$-noncrossing structures, 
despite the fact that they are $D$-finite.

For RNA secondary structures ($2$-noncrossing RNA structures),
certain combinatorial restrictions, for instance minimum arc-length
or stack-size are relatively straightforward to deal with. The
combinatorics and prediction of RNA secondary structures has been 
pioneered by Waterman {\it et al.} in a series of excellent papers
\cite{Penner:93c,Waterman:79a,Waterman:78a,Waterman:94a,Waterman:80}.
He proved for the number of RNA secondary structures of length $n$
(arc-length $\ge 2$), ${\sf T}_2^{[2]}(n)$, the fundamental
recursion
\begin{equation}\label{E:basic}
{\sf T}_2^{[2]}(n)= {\sf T}_2^{[2]}(n-1)+
\sum_{s=0}^{n-3} {\sf T}_2^{[2]}(n-2-s){\sf T}_2^{[2]}(s) \ ,
\end{equation}
where ${\sf T}_2^{[2]}(0)={\sf T}_2^{[2]}(1)={\sf T}_2^{[2]}(2)=1$.
Eq.~(\ref{E:basic}) is an immediate consequence considering
secondary structures as peak-free Motzkin-paths, i.e.,~peak-free
paths with {\it up}, {\it down} and {\it horizontal} steps that stay
in the upper halfplane, starting at the origin and end on the
$x$-axis. The recursion is in particular the key for all asymptotic
results since it immediately implies a functional equation for the
corresponding generating function. This allows the application of
Darboux-type theorems \cite{Schuster:98a,Wong:74}. For the number of
secondary structures with minimum arc-length $\lambda$, ${\sf
T}_2^{[\lambda]}(n)$, it is straightforward to derive
\begin{equation}\label{E:basic2}
{\sf T}_2^{[\lambda]}(n)= {\sf T}_2^{[\lambda]}(n-1)+
\sum_{s=0}^{n-(\lambda+1)}
{\sf T}_2^{[\lambda]}(n-2-s){\sf T}_2^{[\lambda]}(s) \ .
\end{equation}
All asymptotic formulae for secondary structures are of the same type:
a square root. In other words, the asymptotic behavior is determined by an
algebraic branch singularity with the subexponential factor
$n^{-\frac{3}{2}}$.

The situation changes for $k$-noncrossing RNA structures. A
different approach has to be made, since in lack of functional
equations Darboux-type theorems \cite{Wong:74}, cannot be employed.
The idea is to analyze the dominant singularities directly, using
Hankel contours. Singularity analysis has been pioneered 
by P.~Flajolet and A.M.~Odlyzko \cite{Flajolet:07a}. Its basic idea is the 
construction of the ``singular-analogue'' of the Taylor-expansion. 
It can be shown that, under
certain conditions, there exists an approximation, which is locally
of the same order as the original function. The particular, local
approximation allows then to derive the asymptotic form of the
coefficients. In contrast to the subtraction of
singularities-principle \cite{Odlyzko:95a} the only contributions to
the contour integral come from segments close to the singularity. In
our situation all conditions for singularity analysis are satisfied,
since the generating functions involved are $D$-finite 
\cite{Stanley:80,Zeilberger:90}
and $D$-finite functions have an analytic continuation into any
simply-connected domain containing zero. Our approach also works for 
tangled diagrams \cite{Reidys:08Zeil}, which represent the
combinatorial framework for RNA tertiary interactions. 
Our analysis confirms that the particular singularity-type of the 
generating function of $k$-noncrossing RNA
structures depends solely on the crossing number
\cite{Reidys:07asy1,Reidys:lego}. While the location of the singularity
shifts as a function of the arc-length, all subexponential factors remain
the same. Furthermore an interesting feature
is the appearance of logarithms for $k\equiv 1\mod 2$ in the
singular expansion.

\begin{figure}[ht]
\centerline{\epsfig{file=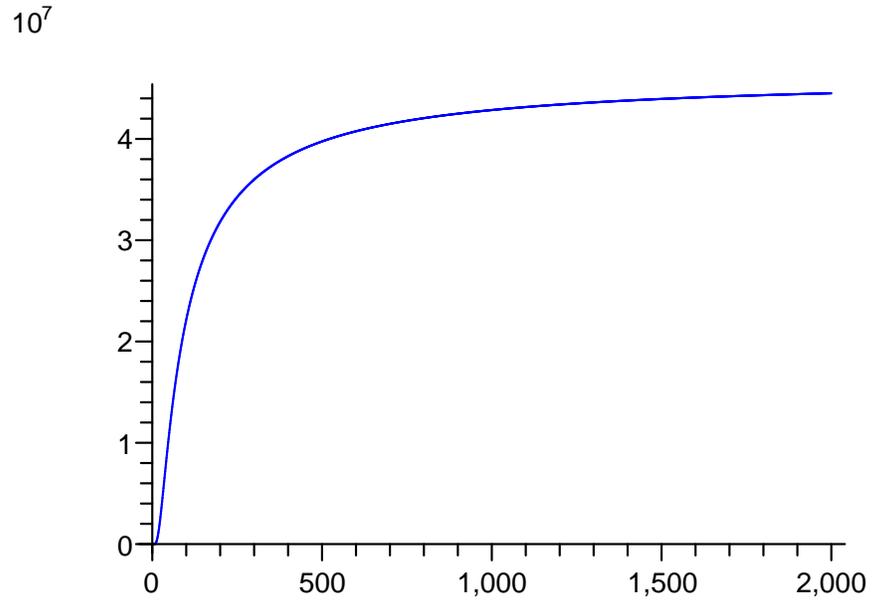,width=1\textwidth}\hskip8pt}
\caption{The ratio $r(n)={\sf T}_{4}^{[4]}/(n^{-21/2}\gamma_4^{-n})$
as a function of $n$. The curve shows that the asymptotic
approximation is valid as $r(n)\sim c_4\approx 4.4509\times
10^7 $.} \label{F:c1}
\end{figure}

Due to biophysical constraints a minimum arc-length of four can be assumed
for minimum free energy RNA structures. The key objective of this paper
is to derive and analyze the generating function for $k$-noncrossing RNA
structures with minimum arc-length $4$, see Table~\ref{T:tab1}.
Based on our results the next step is to compute the subset of
canonical structures, i.e.~the subset of structures with arc-length
$\ge 4$, having no isolated arcs. While it is straightforward to
obtain eq.~(\ref{E:basic2}) from eq.~(\ref{E:basic}) considerable
complication arises, when considering $k$-noncrossing structures with
arc-length $>3$. To understand why, one observes that the number of
ways to place $3$-arcs satisfies a new type of recursion, see
eq.~(\ref{E:recursion}). As a result and in contrast to
$k$-noncrossing structures with minimum arc-length $\lambda\le 3$
the generating function $\sum_{n\ge 0}{\sf T}_k^{[4]}(n)\, z^n$
turns out to be a sum of two power series (Theorem~\ref{T:f4}). The
exponential growth rate can easily be computed via the formula given
in Theorem~\ref{T:asymptotic}, see Table~\ref{T:tab2} and Figure~\ref{F:c1}.

The paper is organized as follows: in Section~\ref{S:preliminaries}
we provide the background on the methods used in this paper. In
Section~\ref{S:functional} we prove a functional equation relating
RNA structures to $k$-noncrossing matchings. We then study the
singularity of the generating function and obtain the asymptotic
formula in Section~\ref{S:asymp}. Finally, in Section~\ref{S:Claim1}
we detail some key ideas instrumental for the proof of Theorem~\ref{T:f4}.

\begin{figure}[ht]
\centerline{%
\epsfig{file=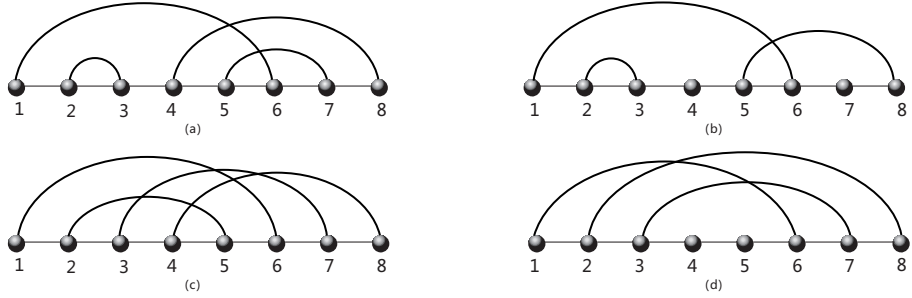,width=0.8\textwidth}\hskip15pt
 }
\caption{\small Basic diagram types: (a) $3$-noncrossing matching
(no isolated points), (b) $3$-noncrossing partial matching (isolated
points $4$ and $7$), (c) $4$-noncrossing RNA structure with
arc-length $\ge 3$, (d) $3$-noncrossing RNA structure with
arc-length $\ge 4$.} \label{F:4-pic}
\end{figure}

\begin{newpage}
\section{Preliminaries}\label{S:preliminaries}
In this Section we provide some background on the generating
functions of $k$-noncrossing matchings \cite{CDDSY,Reidys:08k} and
$k$-noncrossing RNA structures
\cite{Reidys:07pseu,Reidys:07asy1}. We denote the
numbers of $k$-noncrossing matchings and RNA structures with
arc-length $\ge \lambda$ by $f_k(2n)$ and ${\sf
T}_{k}^{[\lambda]}(n)$, respectively. The former corresponds to
$k$-noncrossing diagrams without isolated points and the latter to
$k$-noncrossing diagrams with arc-length $\ge \lambda$. Furthermore,
let ${\sf T}_{k}^{[\lambda]}(n,\ell)$ denote the number of
$k$-noncrossing RNA structures with arc-length $\ge \lambda$ having
exactly $\ell$ isolated points and ${\sf M}_{k}(n)$ denotes the number 
of partial
matchings, or equivalently the number of $k$-noncrossing diagrams
over $[n]$ (i.e.~with isolated points and minimum arc-length $1$).
Pfringsheim's Theorem \cite{Titmarsch:39} guarantees the existence of
a positive real, dominant singularity of $\sum_{n\ge 0}{\sf
M}_k(n)\, z^{n}$ which we denote by $\mu_k$. In order to get some
intuition about the various types of diagrams involved, see
Figure~\ref{F:4-pic}.

\subsection{$k$-noncrossing matchings}
Our main objective is to discuss some basic properties of $f_k(2n)$
and to give an asymptotic formula. Let us recall that a power series $u(x)$
is called $D$-finite over the function field $K(x)$ if 
$\text{\rm dim}\langle u,u', \dots\rangle_{K(x)}<\infty$
\cite{Stanley:80}. The generating function of $k$-noncrossing
matchings satisfies the following identity due to Grabiner {\it
et al.} \cite{Grabiner:93a}
\begin{eqnarray}\label{E:ww0}
\label{E:ww1}
\sum_{n\ge 0} f_{k}(2n)\cdot\frac{z^{2n}}{(2n)!} & = &
\det[I_{i-j}(2z)-I_{i+j}(2z)]|_{i,j=1}^{k-1} \ ,
\end{eqnarray}
where
\begin{equation}\label{E:Bessel}
I_{r}(2z)=\sum_{j \ge 0}\frac{z^{2j+r}}{{j!(r+j)!}}
\end{equation}
denotes the hyperbolic Bessel function of the first kind of order $r$.
Eq.~(\ref{E:ww0}) allows to conclude that
\begin{equation}\label{E:Fk}
F_k(z)=\sum_{n\ge 0}f_k(2n) z^{2n}
\end{equation}
is $D$-finite. Indeed, the hyperbolic Bessel
function \cite{Grabiner:93a} itself is $D$-finite and $D$-finite
functions form an algebra closed under taking Hadamard products
\cite{Stanley:80}. Therefore $D$-finiteness of $F_k(z)$ follows from
eq.~(\ref{E:ww0}). However, beyond the cases $k=2$ and $k=3$,
eq.~(\ref{E:ww0}) does not give directly explicit formulas for
$f_k(2n)$ or $F_k(z)$. For small $k$-values asymptotic formulas can be obtained
using the approximation of the Bessel function
\begin{equation}
I_m(z)= \frac{e^z}{\sqrt{2\pi z}}\left(\sum_{h=0}^{H-1}
\frac{(-1)^h}{h!8^h}\prod_{t=1}^{h}(4m^2-(2t-1)^2)z^{-h}+O(\vert
z\vert^{-H})\right)
\end{equation}
which holds for $-\frac{\pi}{2}< \arg(z)< \frac{\pi}{2}$
\cite{NBS:70}. For arbitrary $k$, systematic analysis of the
determinant $\det[I_{i-j}(2x)-I_{i+j}(2x)]|_{i,j=1}^{k-1}$ by Jin {\it et al.}
\cite{Reidys:08k} shows for arbitrary $k$
\begin{equation}\label{E:theorem}
f_{k}(2n) \, \sim  \, c_k  \, n^{-((k-1)^2+(k-1)/2)}\,
(2(k-1))^{2n},\qquad c_k>0 \ .
\end{equation}
In the following we shall denote the dominant singularity of
$F_k(z)$ by $\rho_k=\frac{1}{2(k-1)}$.

\subsection{$k$-noncrossing RNA structures}
$k$-noncrossing RNA structures are \\
$k$-noncrossing diagrams
satisfying specific arc-length conditions. The latter induce
asymmetries (for instance $1$-arcs are not preserved) which prohibit 
enumeration using Gessel and Zeilberger's 
reflection-principle \cite{GZ} directly (the reflection principle 
implies eq.~(\ref{E:ww0})). For any $k\ge 2$ the numbers of
$k$-noncrossing RNA structures with minimum arc-length $\ge 2$ are
given by \cite{Reidys:07pseu}
\begin{eqnarray}
\label{E:sum}{\sf T}_{k}^{[2]}(n) & = &
\sum_{b=0}^{\lfloor n/2\rfloor}(-1)^{b}{n-b \choose b} {\sf M}_k(n-2b) 
\end{eqnarray}
and we have \cite{Reidys:07asy1}
\begin{eqnarray}
{\sf T}_{k}^{[2]}(n) & \sim &
c_k^{[2]}  \, n^{-((k-1)^2+(k-1)/2)}\,
(\gamma_k^{[2]})^{-n},\qquad c_k^{[2]}>0 \ ,
\end{eqnarray}
where $\gamma_k^{[2]}$ is the unique, solution of minimal modulus of
$\frac{z}{z^2-z+1}=\rho_k$.
For $k$-noncrossing RNA structures with arc-length $\ge 3$ we have
according to \cite{Reidys:07pseu}
\begin{eqnarray}
\label{E:da3}
\forall\, k>2;\qquad
{\sf T}_{k}^{[3]}(n) & = & \sum_{b \le \lfloor \frac{n}{2}\rfloor}
(-1)^b\lambda(n,b)\, {\sf M}_k(n-2b) \ ,
\end{eqnarray}
where $\lambda(n,b)$ denotes the number of way selecting $b$ arcs of
length $\le 2$ over $n$ vertices. The nonexplicit terms
$\lambda(n,b)$ vanish in the functional equation
\begin{align}\label{E:well2}
\sum_{n \ge 0}{\sf T}_{k}^{[3]}(n)\, z^n = \nonumber \\ & 
\frac{1}{1-z+z^2+z^3-z^4}\sum_{n \ge
0}f_k(2n)\left(\frac{z-z^3}{1-z+z^2+z^3-z^4}\right)^{2n} \ .
\end{align}
Singularity analysis based on eq.~(\ref{E:well2}) eventually allows
to derive the asymptotic formula
\begin{eqnarray}
{\sf T}_{k}^{[3]}(n) & \sim &
c_k^{[3]}  \, n^{-((k-1)^2+(k-1)/2)}\,
(\gamma_k^{[3]})^{-n},\qquad c_k^{[3]}>0 \ ,
\end{eqnarray}
where $\gamma_k^{[3]}$ denotes the unique, minimal positive real solution
of $\frac{z-z^3}{1-z+z^2+z^3-z^4}=\rho_k$.

\subsection{Singularity Analysis}
Pfringsheim's Theorem \cite{Titmarsch:39} guarantees that each power
series with positive coefficients has a positive real dominant
singularity. This singularity plays a key role for the asymptotics
of the coefficients. In the proof of Theorem~\ref{T:asymptotic} it
will be important to deduce relations between the coefficients from
functional equations of generating functions. The class of theorems
that deal with such deductions are called transfer-theorems
\cite{Flajolet:07a}. One key ingredient in this framework is a
specific domain in which the functions in question are analytic,
which is ``slightly'' bigger than their respective radius of
convergence. It is tailored for extracting the coefficients via
Cauchy's integral formula. Details on the method can be found in
\cite{Stanley:80,Flajolet:07a}. In case of $D$-finite functions we
have analytic continuation in any simply-connected domain containing
zero \cite{Wasow:87} and all prerequisits of singularity analysis
are met. We use the notation
\begin{equation}\label{E:genau}
\left\{f(z)=O\left(g(z)\right) \
\text{\rm as $z\rightarrow \rho$}\right\}\quad \Longleftrightarrow \quad
\left\{\frac{f(z)}{g(z)} \ \text{\rm is bounded as $z\rightarrow \rho$}\right\}
\end{equation}
The key result used in Theorem~\ref{T:asymptotic} is
\begin{theorem}\label{T:transfer1}\cite{Flajolet:07a}
Let $f(z),g(z)$ be $D$-finite functions with unique dominant singularity
$\rho$ and suppose
\begin{equation}
f(z) = O( g(z)) \ \mbox{ as }\, z\rightarrow \rho \ .
\end{equation}
Then we have
\begin{equation}
[z^n]f(z)= \,C \,\left(1-O(\frac{1}{n})\right)\,  [z^n]g(z)
\end{equation}
where $C$ is a constant and $[z^n]h(z)$ denotes the $n$-th
coefficient of the power series $h(z)$ at $z=0$.
\end{theorem}
\end{newpage}


\section{The generating function}\label{S:functional}


In this Section we compute the generating
function of ${\sf T}_{k}^{[4]}(n)$, the number of $k$-noncrossing
RNA structures with arc-length $\ge 4$. Our first result is a
technical lemma which is instrumental in the proof of
Theorem~\ref{T:f4} below. The proof of the lemma given below is new
and uses integral representations \cite{IRC} instead of dealing with
the combinatorial coefficients directly. 

\begin{lemma}\label{L:laplace} Let $z$ be an
indeterminate over $\mathbb{C}$. Then we have the identity of power
series
\begin{equation}\label{E:1}
\forall \; \vert z\vert < \mu_k;\qquad
 \sum_{n \ge 0}{\sf M}_{k}(n)\, z^n = \left(\frac{1}{1-z}\right)\, \sum_{n \ge
0}f_k(2n)\,\left(\frac{z}{1-z}\right)^{2n} \ .
\end{equation}
\end{lemma}
\begin{proof}
Expressing the combinatorial terms by contour integrals \cite{IRC}
we obtain
\begin{equation}
{{n} \choose {2m}}=\frac{1}{2\pi i}\oint_{|u|=\alpha} (1+u)^n
u^{-2m-1}du \qquad
f_k(2m)=\frac{1}{2 \pi i}\oint_{|v|=\beta}F_k(v)v^{-2m-1}dv
\end{equation}
where $\alpha $, $\beta$ are arbitrary small positive numbers. We derive \\
\begin{eqnarray*}
{\sf M}_{k}(n) &=& \frac{1}{(2\pi i)^2} \sum_m \oint_{|u|=\alpha,
|v|=\beta}(1+u)^n u^{-2m-1}F_k(v)v^{-2m-1}du dv\\
&=& \frac{1}{(2\pi i)^2 } \oint_{|u|=\alpha,|v|=\beta}
(1+u)^n \frac{uv}{(uv)^2-1 } F_k(v) du dv
\end{eqnarray*}
and furthermore
\begin{eqnarray*}
{\sf M}_{k}(n)  &=&\frac{1}{(2\pi i)^2 } \oint_{|v|=\beta} F_k(v)v^{-1}
\left[\oint_{|u|=\alpha} \frac{(1+u)^n
u}{(u+\frac{1}{v})(u-\frac{1}{v})}du\right]
dv \ .\\
\end{eqnarray*}

Since $u=\frac{1}{v}$ and $u=-\frac{1}{v}$ are the only
singularities (poles) enclosed by the particular contour,
eq.~\eqref{E:1} implies

\begin{eqnarray*}
\oint_{|u|=\alpha} \frac{(1+u)^n
u}{(u+\frac{1}{v})(u-\frac{1}{v})}du &=& 2\pi i\left[\frac{(1+u)^n
u}{u-\frac{1}{v}}|_{u=-\frac{1}{v}} +\frac{(1+u)^n
u}{u+\frac{1}{v}}|_{u=\frac{1}{v}}\right]\\
&=&\pi i\, \left(\left[1-\frac{1}{v}\right]^n+\left[1+
\frac{1}{v}\right]^n\right).
\end{eqnarray*}
Therefore, for $\vert z\vert <\mu_k$
\begin{eqnarray*}
\sum_{n \ge 0}{\sf M}_{k}(n)z^n &=& \frac{1}{4\pi i}\sum_{n \ge 0}
\oint_{|v|=\beta} F_k(v)v^{-1}\left(\left[1-\frac{1}{v}\right]^n+\left[1+
\frac{1}{v}\right]^n\right)z^n dv\\
&=&\frac{1}{4\pi i}\oint_{|v|=\beta} F_k(v)\frac{1}{v-(v-1)z}dv +
\frac{1}{4\pi i}\oint_{|v|=\beta} F_k(v)\frac{1}{v-(v+1)z}dv \ .
\end{eqnarray*}
The first integrand has its unique pole at $v=-\frac{z}{1-z}$ and
the second at $v=\frac{z}{1-z}$, respectively:
$$
\frac{1}{v-(v-1)z}=\frac{1}{v+\frac{z}{1-z}}\, \frac{1}{1-z} \quad
\text{\rm and}\quad
\frac{1}{v-(v+1)z} =\frac{1}{v-\frac{z}{1-z}}\, \frac{1}{1-z} \ .
$$
In view of $F_k(z)=F_k(-z)$ we derive
\begin{eqnarray*}
\sum_{n \ge 0}{\sf M}_{k}(n)z^n =
\frac{1}{1-z}\left[\frac{1}{2}F_k\left(-\frac{z}{1-z}\right)+
\frac{1}{2}F_k\left(\frac{z}{1-z}\right)\right] =
\frac{1}{1-z}F_k\left(\frac{z}{1-z}\right) \ ,
\end{eqnarray*}
whence the lemma.
\end{proof}

Before we state the main result of this section, let us introduce
some notation. We set

\begin{eqnarray}
u(z)   & = & \sqrt{1+4z-4z^2-6z^3+4z^4+z^6} \\
\label{E:oben}
f_j(z) & = & -\frac{-2z^2+z^3-1+(-1)^{j}\,u(z)}{2(1-2z-z^2+z^4)} \ .
\end{eqnarray}
Note that $f_j(z)$ is an algebraic function over the function field $K(z)$,
i.e.~there exists a polynomial with coefficients being polynomials in $z$ for
which $f_j(z)$ is a root. This fact will be important when computing the 
subexponential factors of the asymptotic formula for ${\sf T}_k^{[4]}(n)$.
\begin{theorem}\label{T:f4}
Let $k$ be a positive integer, $k>3$ and $f_{1}(z)$ and $f_2(z)$ be given by
eq.(\ref{E:oben}). Then we have the functional equation
\begin{eqnarray*}\label{E:f4}
\sum_{n \ge 0}{\sf T}_{k}^{[4]}(n)\, z^n &  = &
\frac{F_{1}(-z^2)}{1-zf_{1}(-z^2)}
\sum_{n \ge
0}f_k(2n)\left(\frac{z\, f_{1}(-z^2)}{1-zf_{1}(-z^2)}\right)^{2n}+\\
& & \frac{F_{2}(-z^2)}{1-zf_{2}(-z^2)}\sum_{n \ge
0}f_k(2n)\left(\frac{z\,f_{2}(-z^2)}{1-zf_{2}(-z^2)}\right)^{2n} \ .
\end{eqnarray*}
\end{theorem}
\begin{proof}
{\it Claim $1$.} Let $\lambda(n,b)$ denote the number of ways to place $b$
arcs of length $\le 4$ over $[n]$. Then we have
\begin{eqnarray}\label{E:KK}
{\sf T}_{k}^{[4]}(n) & = & \sum_{b \le \lfloor \frac{n}{2}\rfloor}
(-1)^b\, \lambda(n,b)\, {\sf M}_k(n-2b)
\end{eqnarray}

and $\lambda(n,b)$ satisfies the recursion

\begin{equation}\label{E:recursion}
\begin{split}
&\lambda(n+2b,b)= \\ &\lambda(n+2b-1,b)+\lambda(n+2b-4,b-2)+
\lambda(n+2b-5,b-2)+\lambda(n+2b-6,b-3)\\
&+\sum_{i=1}^{b}{[\lambda(n+2b-2i,b-i)+2\lambda(n+2b-2i-1,b-i)+
\lambda(n+2b-2i-2,b-i)]}\\
&-\lambda(n+2b-3,b-1)\ ,
\end{split}
\end{equation}
where $\lambda(n,0)=1$, $\lambda(n,1)=3n-6$
and $n \geq 2b$. The proof of Claim $1$ is analogous to the proof of
Theorem~$5$ in \cite{Reidys:07pseu}. In order to keep the paper
selfcontained we present it in Section~\ref{S:Claim1}. \\
The idea is
now to relate $\sum_{n \ge 0}{\sf T}_{k}^{[4]}(n)\, z^n$ to the
power series $\sum_{n\ge 0}{\sf M}_k(n)\,z^n$. For this purpose we
compute
\begin{align*}
\sum_{n \ge 0}{\sf T}_{k}^{[4]}(n)z^n&=\sum_{n \ge 0}\sum_{2b\le n}
(-1)^b\lambda(n,b)\sum_{m=2b}^{n}{n-2b \choose
m-2b}f_k(m-2b,0)\, z^n\\
&=\sum_{b \ge 0}(-1)^b z^{2b}\sum_{n \ge
2b}\lambda(n,b)\, {\sf M}_{k}(n-2b)z^{n-2b}\\
&=\sum_{b \ge 0}(-1)^b z^{2b}\sum_{n \ge 0}\lambda(n+2b,b)\, {\sf M}_{k}(n)
\,z^n \ .
\end{align*}
Interchanging the summations w.r.t.~$b$ and $n$ we arrive at
\begin{equation}\label{E:21}
\sum_{n \ge 0}{\sf T}_{k}^{[4]}(n)z^n = \sum_{n \ge 0}\left[\sum_{b\ge
0}(-1)^b z^{2b}\lambda(n+2b,b)\right]{\sf M}_{k}(n)\, z^n \ .
\end{equation}
Now we use the recursion formula for $\lambda(n,b)$.
Let
\begin{equation}\label{E:varphi}
\varphi_{n}(z)=\sum_{b \ge 0}\lambda(n+2b,b)z^b \ .
\end{equation}
Multiplying in eq.~(\ref{E:recursion}) with
$z^b$ and taking the summation over all $b$ ranging from $0$ to $\lfloor
n/2\rfloor$ implies for $\varphi_{n}(z)$,
$n=1,2\ldots$
\begin{equation}\label{E:22}
\left(1-z^2-z^3-\frac{z}{1-z}\right)\varphi_n(z)=\left(z^2+\frac{z^2+1}{1-z}
\right)\varphi_{n-1}(z)+\left(\frac{z}{1-z}\right)\varphi_{n-2}(z) \ .
\end{equation}
We make the Ansatz
\begin{equation}\label{E:ansatz}
f(x,y)=\sum_{n \ge 0}\sum_{j \le
\frac{n}{2}}\lambda(n,j)x^j\ \frac{y^n}{n!}=\sum_{n \ge
0}\varphi_n(x)\ \frac{y^n}{n!} \ .
\end{equation}
Multiplying in eq.~(\ref{E:22}) with
$\frac{y^n}{n!}$ and taking the summation over all $n \ge 0$ leads to
the partial differential equation
\begin{equation}\label{E:diff}
\left(1-x^2-x^3-\frac{x}{1-x}\right)\frac{\partial^2f(x,y)}{\partial
y^2}=\left(x^2+\frac{x^2+1}{1-x}\right)\frac{\partial f(x,y)}{\partial
y}+\left(\frac{x}{1-x}\right)f(x,y) \ .
\end{equation}
The general solution of eq.~(\ref{E:diff}) can be computed by MAPLE and
is given by
\begin{align*}
f(x,y)&= F_{1}(x)\exp(f_1(x)\cdot y)+F_2(x)\exp(f_2(x)\cdot y) \\
&=\sum_{n \ge
0}\left[F_1(x)\,f_1(x)^n+F_2(x)\,f_2(x)^n\right]\frac{y^n}{n!} \ ,
\end{align*}
where $F_1(x)$, $F_2(x)$ are arbitrary functions and
\begin{equation}
f_1(x)=\frac{2x^2-x^3+1+u(x)}{2(1-2x-x^2+x^4)}, \quad
f_2(x)=\frac{2x^2-x^3+1-u(x)}{2(1-2x-x^2+x^4)} \ .
\end{equation}
By definition we have $f(x,y)=\sum_{n \ge 0}\varphi_n(x)\cdot\frac{y^n}{n!}$
and
\begin{equation}\label{E:explicit}
\varphi_n(x)=F_1(x)(f_1(x))^n+F_2(x)(f_2(x))^n \ .
\end{equation}
In order to solve eq.~(\ref{E:explicit}) it remains to compute $F_1(x)$ and
$F_2(x)$. The key information lies in the
initial conditions for $f(x,y)$ and $\varphi_{n}(x)$.
Explicitly we have $f(x,0)=1$ and $\varphi_{1}(x)=
\lambda(1,0)\, x^0=1$, which implies
\begin{eqnarray*}
F_1(x)+F_2(x) & = & 1\\
F_1(x)f_1(x)+F_2(x)f_2(x) & = & 1 \ .
\end{eqnarray*}
Accordingly we obtain
\begin{equation}
F_1(x)=\frac{f_2(x)-1}{f_2(x)-f_1(x)}\quad \text{\rm and}\quad
F_2(x)=\frac{f_1(x)-1}{f_1(x)-f_2(x)} \ .
\end{equation}
In view of $\varphi_{n}(-z^2)  =  \sum_{b \ge 0}\lambda(n+2b,b)(-1)^bz^{2b}$
we can express $\sum_{n \ge 0}{\sf T}_{k}^{[4]}(n)z^n$ as follows:
\begin{eqnarray*}
\sum_{n \ge 0}{\sf T}_{k}^{[4]}(n)z^n & = &
\sum_{n \ge 0}\varphi_{n}(-z^2) \, {\sf M}_{k}(n)
\,z^n  \\
 & = &
F_1(-z^2)\, \sum_{n \ge 0} {\sf M}_{k}(n) \, \left(f_1(-z^2)z\right)^n \, +\,
F_2(-z^2)\, \sum_{n \ge 0} {\sf M}_{k}(n) \, \left(f_2(-z^2)z\right)^n \ .
\end{eqnarray*}
Now we use Lemma~\ref{L:laplace}:
$$
\sum_{n \ge 0}{\sf M}_{k}(n)\, z^n = \left(\frac{1}{1-z}\right)\, \sum_{n \ge
0}f_k(2n)\,\left(\frac{z}{1-z}\right)^{2n} \ ,
$$
which allows to express $\sum_{n \ge 0}{\sf T}_{k}^{[4]}(n)z^n$ via
$\sum_{n\ge 0}f_k(2n)\,z^{2n}$
\begin{align*}\label{E:laplace}
\sum_{n \ge 0}{\sf T}_{k}^{[4]}(n)\, z^n
=\frac{F_{1}(-z^2)}{1-zf_{1}(-z^2)} \sum_{n \ge
0}f_k(2n)\left(\frac{zf_{1}(-z^2)}{1-zf_{1}(-z^2)}\right)^{2n}+\\
\frac{F_{2}(-z^2)}{1-zf_{2}(-z^2)} \sum_{n \ge
0}f_k(2n)\left(\frac{zf_{2}(-z^2)}{1-zf_{2}(-z^2)}\right)^{2n}\  .
\end{align*}
\end{proof}


\section{Asymptotics of RNA pseudoknot structures with arc-length $\ge 4$}
\label{S:asymp}

We set
\begin{eqnarray}\label{E:t1}
\vartheta_1(z) & = & \frac{z\, f_{1}(-z^2)}{1-zf_{1}(-z^2)} \\
\label{E:t2}
\vartheta_2(z) & = & \frac{z\, f_{2}(-z^2)}{1-zf_{2}(-z^2)} \ .
\end{eqnarray}
Note that $\vartheta_1(z)$ and $\vartheta_2(z)$ are algebraic functions
over the function field $K(z)$.

\begin{theorem}\label{T:asymptotic}
Let $k>3$ be a positive integer and $\rho_k,\gamma_k$ denote the
positive real singularities of $F_k(z)=\sum_{n\geq 0}f_k(2n)z^{2n}$
and $\sum_{n\geq 0}{\sf T}_{k}^{[4]}(n)\, z^n$, respectively. Then
the number of $k$-noncrossing RNA structures with arc-length $\geq
4$ is for $k\le 9$ asymptotically given by
\begin{equation}
 {\sf T}_{k}^{[4]}(n) \sim c_k\, n^{-((k-1)^2+(k-1)/2)}\,
\left(\gamma_k^{-1}\right)^n \ ,
\end{equation}
where $\gamma_k$ is the unique positive, real solution of the equation
$\vartheta_1(z)=\rho_k$.
\end{theorem}
\begin{proof}
According to Theorem~\ref{T:f4} we have the functional equation
\begin{eqnarray*}
\sum_{n \ge 0}{\sf T}_{k}^{[4]}(n)\, z^n &  = &
\frac{F_{1}(-z^2)}{1-zf_{1}(-z^2)}\, \underbrace{\sum_{n \ge
0}f_k(2n)\left(\frac{z\, f_{1}(-z^2)}{1-zf_{1}(-z^2)}\right)^{2n}}_{
F_k(\vartheta_1(z))}+\\
& & \frac{F_{2}(-z^2)}{1-zf_{2}(-z^2)}\, \underbrace{\sum_{n \ge
0}f_k(2n)\left(\frac{z\,f_{2}(-z^2)}{1-zf_{2}(-z^2)}\right)^{2n}}_{
F_k(\vartheta_2(z))}
\ .
\end{eqnarray*}
We consider the functions $\vartheta_1(z)$, $\vartheta_2(z)$ given
by eq.~(\ref{E:t1}) and eq.~(\ref{E:t2}). The mappings $x\mapsto
\vartheta_1(x)$ and $x\mapsto \vartheta_2(x)$ are strictly monotone and
$\vartheta_1(x)> \vartheta_2(x)$ for $\vartheta_1(x)\in ]0,\frac{1}{5}]$.
Furthermore we have $\rho_k < \rho_4=\frac{1}{6}$, for $k>4$.
We can conclude from this that the real, positive dominant
singularity, $\gamma_k$, of $\sum_{n \ge 0}{\sf T}_{k}^{[4]}(n)\, z^n$,
whose existence is guaranteed by Pfringsheim's Theorem \cite{Titmarsch:39},
satisfies
\begin{equation}
\vartheta_1(\gamma_k)=\rho_k \ .
\end{equation}
Being a determinant of Bessel functions \cite{Grabiner:93a},
$F_k(z)$ is $D$-finite. Moreover $\vartheta_1(z)$ and $\vartheta_2(z)$ are
algebraic over $K(z)$, analytic for $\vert z\vert <\delta$, where
$\gamma_k<\delta$ and satisfy $\vartheta_1(0)=\vartheta_2(0)=0$.
Therefore the composition $F_k(\vartheta_i(z))$, $i=1,2$, is
$D$-finite \cite{Stanley:80} and $F_k(\vartheta_1(z))$ and
$F_k(\vartheta_2(z))$ have singular expansions, respectively. We
further observe that neither $\frac{F_{1}(-z^2)}{1-zf_{1}(-z^2)}$
nor $\frac{F_{2}(-z^2)}{1-zf_{2}(-z^2)}$ have a singularity $\zeta$
with $\vert \zeta\vert \le \gamma_k$. Hence if $\zeta$ is a dominant
singularity of $\sum_{n}{\sf T}_k^{[4]}(n)\,z^n$ then it is
necessarily a singularity of $F_k(\vartheta_1(z))$ or
$F_k(\vartheta_2(z))$. As for singularities of $F_k(\vartheta_1(z))$
and $F_k(\vartheta_2(z))$, we consider for $k\le 9$ the ODE
satisfied by $F_k(z)$:
\begin{equation}\label{E:JK}
q_{0,k}(z)\frac
{d^e}{dz^e}F_k(z)+q_{1,k}(z)\frac{d^{e-1}}{dz^{e-1}}F_k(z)+\dots +
q_{e,k}(z)F_k(z)=0  \  ,
\end{equation}
where $q_{j,k}(z)$ are polynomials. The key point is now that any
dominant singularity of $F_k(z)$ is contained in the set of roots of
$q_{0,k}(z)$ \cite{Stanley:80}. Computing the ODEs for $4\le k\le 9$
we can therefore conclude that $F_k(z)$ has only the two dominant
singularities $\rho_k$ and $-\rho_k$. Let $S=\{\zeta\mid
\vartheta_1(\zeta) = \rho_k \;\text{\rm or}\; \vartheta_2(\zeta) =
-\rho_k\}$. Then $\gamma_k$ is the unique $S$-element of minimal
modulus. We can draw two conclusions: first
\begin{equation}\label{E:mm}
[z^n]\, {\sf T}_k^{[4]}(z)\,\sim  \, c_k
\, [z^n]\, F_k(\vartheta_1(z))\quad \text{\rm for some $c_k>0$}
\end{equation}
and secondly, $\gamma_k$ is the unique dominant singularity of
$\sum_{n}{\sf T}_k^{[4]}(n)\,z^n$. In view of eq.~(\ref{E:mm}) it
thus remains to analyze the subexponential factors of the singular
expansion of $F_k(\vartheta_1(z))$ at $z=\gamma_k$. Since
$\vartheta_1(z)$ is regular at $\gamma_k$ we are given the
supercritical case of singularity analysis \cite{Flajolet:07a}. In
the supercritical case the subexponential factors of the compositum,
$F_k(\vartheta_1(z))$ coincide with those of the outer function,
$F_k(z)$. According to \cite{Reidys:08k} we have for arbitrary $k$
\begin{equation}
f_k(2n)\sim n^{-((k-1)^2 +\frac{k-1}{2})}\,
\left({\rho_k}^{-1}\right)^{2n}
\end{equation}
and therefore the subexponential factors of 
$F_k(z)=\sum_{n\ge 0}f_k(2n)z^{2n}$ coincide with those of $F_k(\vartheta_1(z))$,
i.e.~we have
\begin{equation}
{\sf T}_k^{[4]}(n)
\sim c_k \, n^{-((k-1)^2 +\frac{k-1}{2})}\,\left({\gamma_k}^{-1}\right)^n 
\end{equation}
proving the theorem.
\end{proof}

\section{Proof of Claim $1$}\label{S:Claim1}
We recall that the numbers of $k$-noncrossing matchings and RNA
structures with arc-length $\ge \lambda$ are denoted by $f_k(2n)$
and ${\sf T}_{k}^{[\lambda]}(n)$, respectively. Furthermore,
${\sf T}_{k}^{[\lambda]}(n,\ell)$ denotes the number of
$k$-noncrossing RNA structures with arc-length $\ge \lambda$ having
exactly $\ell$ isolated points, and let $f_k(m,\ell)$ denote the number of
$k$-noncrossing diagrams with $\ell$ isolated points over $m$
vertices. Let $\mathscr{G}_{n,k}(\ell,j_1,j_2,j_3)$ be the set of
all $k$-noncrossing diagrams having exactly $\ell$ isolated points
and exactly $j_1$ $1$-arcs, $j_2$ $2$-arcs and $j_3$ $3$-arcs. We
set $G_k(n,\ell,j_1,j_2,j_3)=|\mathscr{G}_{n,k}(\ell,j_1,j_2,j_3)|$.
In particular, we have $G_k(n,\ell,0,0,0)={\sf
T}_{k}^{[4]}(n,\ell)$. We observe that Claim $1$ is implied (taking
the sum over all $\ell$) by
\begin{eqnarray}\label{E:ww}
{\sf T}_{k}^{[4]}(n,\ell) & = & \sum_{b \le \lfloor \frac{n}{2}\rfloor}
(-1)^b\lambda(n,b)\, f_{k}(n-2b,\ell) \ ,
\end{eqnarray}
where $\lambda(n,b)$ satisfies the recursion
\begin{equation}\label{E:rec2}
\begin{split}
\lambda(n,b) &=\lambda(n-1,b)+\lambda(n-4,b-2)+
\lambda(n-5,b-2)+\lambda(n-6,b-3) \\
&+\sum_{i=1}^{b}{[\lambda(n-2i,b-i)+2\lambda(n-2i-1,b-i)+
\lambda(n-2i-2,b-i)]}\\
&-\lambda(n-3,b-1)
\end{split}
\end{equation}
with the initial conditions $\lambda(n,0)=1$, $\lambda(n,1)=3n-6$ and $n
\geq 2b$.\\
We shall proceed by proving eq.~(\ref{E:ww}). For this purpose, let
$\lambda(n,b_1,b_2,b_3)$ denote the number
of ways to select exactly $b_1$ $1$-arcs, $b_2$ $2$-arcs and $b_3$ $3$-arcs
over ${1,\ldots,n}$ vertices.\\
{\it Claim A.}
\begin{equation}
\sum_{j_1\geq b_1,j_2 \geq b_2,j_3 \geq b_3} {j_1 \choose b_1} {j_2
\choose b_2} {j_3 \choose b_3}
G_k(n,\ell,j_1,j_2,j_3)=\lambda(n,b_1,b_1,b_3)f_k(n-2(b_1+b_2+b_3),\ell).
\end{equation}
The idea is to construct a family $\mathcal{F}$ of
$\mathscr{G}_{n,k}$-diagrams, having $\ell$ isolated points and {
at least} $b_1$ $1$-arcs, $b_2$ $2$-arcs and $b_3$ $3$-arcs,
respectively. We then express $|\mathcal {F}|$ via the numbers
$G_k(n,\ell,j_1,j_2,j_3)$. We select (a) $b_1$ $1$-arcs and $b_2$
$2$-arcs and $b_3$ $3$-arcs and (b) an arbitrary $k$-noncrossing
diagram over the remaining $n-2(b_1+b_2+b_3)$ vertices with exactly
$\ell$ isolated points. Let $\mathcal {F}$ be the family of diagrams
obtained in this way. It is straightforward to show that
$\lambda(n,b_1,b_2,b_3)$ satisfies the recursion:
\begin{align*}\label{E:lambda}
&\lambda(n,b_1,b_2,b_3)=\\  
&\quad\lambda(n-1,b_1,b_2,b_3)+\lambda(n-2,b_1-1,b_2,b_3)+
\lambda(n-4,b_1-1,b_2,b_3-1)\nonumber \\
&+\lambda(n-5,b_1,b_2,b_3-2)+
\lambda(n-6,b_1,b_2,b_3-3)-\lambda(n-3,b_1,b_2-1,b_3)\nonumber \\
&+\sum_{i=1}^{b}{[2\lambda(n-2i-1,b_1,b_2-1,b_3-(i-1))+
\lambda(n-2i-2,b_1,b_2,b_3-i)]}\nonumber \\
&+\sum_{i=2}^{b}{[\lambda(n-2i,b_1,b_2-2,b_3-(i-2))]} \nonumber
\end{align*}
with the initial conditions $\lambda(n,0,0,0)=1$,
$\lambda(n,1,0,0)=n-1$, $\lambda(n,0,1,0)=n-2$,
$\lambda(n,0,0,1)=n-3$, $n\geq 2b$.\\
Clearly, each element $\theta \in \mathcal {F}$ is
contained in $\mathscr{G}_{n,k}(\ell,j_1,j_2,j_3)$ for some $j_1\geq b_1$ and
$j_2\geq b_2$ and $j_3\geq b_3$. Indeed, any $1$-arc or $2$-arc or $3$-arc can
only cross at most two other arcs. Therefore $1$-arcs and $2$-arcs
and $3$-arcs cannot be contained in a set of more than $3$-mutually
crossing arcs. As a result, for $k>3$ the construction generates
$k$-noncrossing diagrams. Clearly, $\theta$ has exactly $\ell$ isolated
vertices and in step (b) we potentially derive additional $1$-arcs
and $2$-arcs and $3$-arcs, whence $j_1\geq b_1$
and $j_2 \geq b_2$ and $j_3 \geq b_3$, respectively. Next we observe that
we have by construction
\begin{eqnarray*}
|\mathcal {F}|=\lambda(n,b_1,b_2,b_3)\, f_k(n-2(b_1+b_2+b_3),\ell) \ .
\end{eqnarray*}
Since any of the $k$-noncrossing diagrams over $n-2(b_1+b_2+b_3)$ vertices can
generate additional $1$-arcs or $2$-arcs or $3$-arcs, we consider
\begin{equation*}
\mathcal {F}(j_1,j_2,j_3)= \{ \theta \in \mathcal {F}\mid
\theta \mbox{ has exactly $j_1$ $1$-arcs, $j_2$ $2$-arcs and  $j_3$ $3$-arcs}
\}.\\
\end{equation*}
Obviously, we then have the partition $\mathcal {F}=\dot\cup_ {j_1\geq
b_1,j_2\geq b_2,j_3\geq b_3} \mathcal{F}(j_1,j_2,j_3)$. Suppose
$\theta \in \mathcal {F}(j_1,j_2,j_3)$, then $\theta \in
\mathscr{G}_{n,k}(\ell,j_1,j_2,j_3)$ and furthermore
$\theta$ occurs with multiplicity  ${j_1 \choose b_1}$ ${j_2 \choose
b_2}$ ${j_3 \choose b_3}$ in $\mathcal {F}$ since by construction
any $b_1$-element subset of the $j_1$ $1$-arcs and $b_2$-element
subset of the $j_2$ $2$-arcs and $b_3$-element subset of the $j_3$
$3$-arcs is counted respectively in $\mathcal {F}$. Therefore we have
\begin{equation}
|\mathcal {F}(j_1,j_2,j_3)| ={j_1 \choose b_1} {j_2 \choose b_2}
{j_3 \choose b_3} G_k(n,\ell,j_1,j_2,j_3)
\end{equation}
and
\begin{align*}
\sum_{j_1\geq b_1,j_2\geq b_2,j_3\geq b_3} |\mathcal {F}(j_1,j_2,j_3)|
&=\lambda(n,b_1,b_2,b_3)f_k(n-2(b_1+b_2+b_3),\ell)
\end{align*}
proving Claim $A$. We next set
$$
F_k(x,y,z)=\sum_{j_1\geq 0}\sum_{j_2 \geq 0}\sum_{j_3
\geq 0}G_k(n,\ell,j_1,j_2,j_3)x^{j_1}y^{j_2}z^{j_3} \ .
$$
Taking derivatives we obtain
\begin{align*}
&\frac{1}{b_1 !}\frac{1}{b_2 !}\frac{1}{b_3 !}F_k^{b_1,b_2,b_3}(1)\\
=&\sum_{j_1\geq b_1,j_2 \geq b_2,j_3 \geq b_3} {j_1 \choose b_1}
{j_2 \choose b_2} {j_3 \choose b_3}
G_k(n,\ell,j_1,j_2,j_3)1^{j_1-b_1}1^{j_2-b_2}1^{j_3-b_3}
\end{align*}
and accordingly
\begin{align*}
& \sum_{j_1\geq 0,j_2 \geq 0,j_3 \geq
0}G_k(n,\ell,j_1,j_2,j_3)x^{j_1}y^{j_2}z^{j_3}\\
& = \sum_{b_1\geq 0,b_2 \geq 0,b_3 \geq 0}\left[\sum_{j_1\geq
b_1,j_2 \geq b_2,j_3 \geq b_3}{j_1 \choose b_1} {j_2 \choose b_2}
{j_3 \choose b_3}
G_k(n,\ell,j_1,j_2,j_3)\right]\\&(x-1)^{b_1}(y-1)^{b_2}(z-1)^{b_3} \\
& = \sum_{b_1\geq 0,b_2 \geq 0,b_3 \geq 0}\lambda(n,b_1,b_2,b_3)
f_k(n-2(b_1+b_2+b_3),\ell)(x-1)^{b_1}(y-1)^{b_2}(z-1)^{b_3} \ .
\end{align*}
By construction $G(n,\ell,0,0,0)$ is the constant term of the
$F_k(x,y,z)$. That is, the number of $k$-noncrossing RNA structures
with $\ell$ isolated vertices and no $1$-arcs, $2$-arcs and $3$-arcs
is given by
\begin{equation}
G(n,\ell,0,0,0)= \sum_{b_1\geq 0,b_2 \geq 0,b_3 \geq
0}(-1)^{b_1+b_2+b_3}\lambda(n,b_1,b_2,b_3)f_k(n-2(b_1+b_2+b_3),\ell) \ .
\end{equation}
We take the sum over all $\ell$ and derive
\begin{align}
&{\sf T}_{k}^{[4]}(n)= \\
&\sum_{b_1 \geq 0,b_2 \geq 0,b_3 \geq 0}^{\lfloor
\frac{n}{2}\rfloor} (-1)^{b_1+b_2+b_3}\lambda(n,b_1,b_2,b_3)
\left[\sum_{\ell=0}^{n-2(b_1+b_2+b_3)}f_{k}(n-2(b_1+b_2+b_3),\ell)\right]
\ .\nonumber
\end{align}
Setting
\begin{equation*}
\lambda(n,b)=\sum_{b_1+b_2+b_3=b}\lambda(n,b_1,b_2,b_3)
\end{equation*}
we conclude first
\begin{eqnarray*}
{\sf T}_{k}^{[4]}(n)
& = & \sum_{b \le \lfloor \frac{n}{2}\rfloor}
(-1)^b\lambda(n,b)\,{\sf M}_k(n-2b)
\end{eqnarray*}
and second eq.~(\ref{E:rec2}), completing the proof of Claim
$1$.\\


{\bf Acknowledgments.}
We are grateful to Fenix W.D.~Huang, Emma Y.~Jin, Jing Qin and Rita
R.~Wang for their help. This work was supported by the 973 Project,
the PCSIRT Project of the Ministry of Education, the Ministry of
Science and Technology, and the National Science Foundation of
China.
\bibliographystyle{amsplain}

\newpage
\begin{table}[!h]
\tabcolsep 0pt
\begin{center}
\def\temptablewidth{1\textwidth}
{\rule{\temptablewidth}{1pt}}
\begin{tabular*}{\temptablewidth}{@{\extracolsep{\fill}}llllllllllllllllll}
    & $1$ & $2$ & $3$ & $4$ & $5$ & $6$ & $7$ & $8$ & $9$ & $10$ & $11$ 
& $12$ & $13$ & $14$ & $15$ \\   \hline
${\sf T}_k^{[4]}(n)$ & $1$ & $1$ & $1$ & $1$ & $2$
             & $5$ & $15$ & $51$ & $179$ & $647$ & $2397$
             & $9081$ & $35181$ & $139307$ & $563218$ \\  \hline
\end{tabular*}
{\rule{\temptablewidth}{1pt}}
\end{center}
\caption{The first $15$ numbers of $4$-noncrossing RNA structures
with arc-length $\geq 4$} \label{T:tab1}\vspace*{-12pt}
\end{table}

\centerline{} \vspace{1.5cm}
\begin{table}[!h]
\tabcolsep 0pt
\begin{center}
\def\temptablewidth{1\textwidth}
{\rule{\temptablewidth}{1pt}}
\begin{tabular*}{\temptablewidth}{@{\extracolsep{\fill}}ccccccc}
 $k$   & $4$ &$5$ &$6$ &$7$ &
$8$  \\   \hline\\ $\gamma^{-1}_k$& $6.52900$ & $8.64830$ &
$10.71759$ & $12.76349$ & $14.79631$
\\ \\
${\sf T}_k^{[4]}(n)$ &$c_4 n^{-\frac{21}{2}}(\gamma^{-1}_4)^n$&$c_5
n^{-18}(\gamma^{-1}_5)^n$&$c_6
n^{-\frac{55}{2}}(\gamma^{-1}_6)^n$&$c_7
n^{-39}(\gamma^{-1}_7)^n$&$c_8 n^{-\frac{105}{2}}(\gamma^{-1}_8)^n$
\\ \\   \hline
\end{tabular*}
{\rule{\temptablewidth}{1pt}}
\end{center}
\caption{Exponential growth rates and asymptotic formulas for
$k$-noncrossing RNA structures with minimum arc-length $\ge 4$.}
\vspace*{-11pt} \label{T:tab2}
\end{table}

\vspace{2.5cm}

\centerline{}

\centerline{}

\end{document}